\documentclass[10pt,reqno]{amsart}
\usepackage{amsmath}
\usepackage{amsfonts}
\usepackage{mathrsfs}
\usepackage{amssymb}
\usepackage{amsthm}
\usepackage{enumitem}
\usepackage{bbm}
\usepackage{times}
\usepackage{tabularx}
\usepackage{amsbsy}
\usepackage{mathtools}
\usepackage{tikz}
\usetikzlibrary{arrows,decorations.markings}
\usetikzlibrary{matrix}
\usetikzlibrary{graphs}
\usetikzlibrary{backgrounds}
\usepackage{setspace}     \spacing{1}
\usepackage{color}

\usepackage[colorlinks=true,linkcolor=blue,citecolor=blue]{hyperref}

\theoremstyle{Theorem}
\newtheorem{theorem}{Theorem} [section]

\newtheorem{proposition}[theorem]{Proposition}

\theoremstyle{definition}
\newtheorem{definition}[theorem]{Definition}

\theoremstyle{remark}

\newlist{enumlemma}{enumerate}{3}
\setlist[enumlemma]{label*={(\alph*)}, ref= {(\alph*)} }

\usepackage{marginnote}

\newcommand{\diff}{\mathrm{Diff}}

\renewcommand{\epsilon}{\varepsilon}

\DeclareMathOperator{\Diff}{Diff}

\DeclareMathOperator{\Homeo}{Homeo}

\newcommand{\ci}{C^{\infty}}

\newcommand{\R}{\mathbb {R}}
\newcommand{\Q}{\mathbb {Q}}
\newcommand{\Z}{\mathbb {Z}}
\newcommand{\N}{\mathbb {N}}
\newcommand{\T}{\mathbb {T}}

\newcommand{\Sl}{\mathrm{SL}}
\newcommand{\PSl}{\mathrm{PSL}}

\def\wl{{l}}

\title{Recent progress in the Zimmer program}

\author[D.~Fisher]{David Fisher}
\address{Indiana University, Bloomington, Bloomington, IN 47401, USA}
\email{fisherdm@indiana.edu}

\thanks{DF was partially supported by NSF Grant  DMS-1607041.}

\long\def\symbolfootnote[#1]#2{\begingroup\def\thefootnote{\fnsymbol{footnote}}
\footnote[#1]{#2}\endgroup}

\begin{document}
\maketitle

\begin{abstract}
This paper can be viewed as a sequel to the author's long survey on the Zimmer program \cite{F11}
published in 2011.  The sequel focuses on recent rapid progress on certain aspects of the program
particularly concerning rigidity of Anosov actions and Zimmer's conjecture that there are no
actions in low dimensions.  Some emphasis is put on the surprising connections between these
two different sets of developments and also on the key connections and ideas for future research that
arise from these works taken together.
\end{abstract}





\section{Introduction}
\label{sec:intro}

One could attempt to write an afterward to this volume that addressed the full scope and impact of Zimmer's work on mathematics.
The result would either be terse and unreadable or another volume longer than the curent one.  While Zimmer's contributions and
ideas have fostered research in many directions, there is an essential unity to his perspective.  So rather than attempting to
broadly survey all the impacts, I chose here to focus on a ``test case", where Zimmer's contributions are particularly important
and where there has been dramatic recent progress: the {\em Zimmer
program}.

The {\em Zimmer program} aims to classify actions of higher rank Lie groups and their lattices on compact manifolds.  The program was initiated by Zimmer in a series of papers in the early 1980's and framed explicitly in \cite{MR682830, MR934329}.  In the last 35 years, numerous authors have
made deep and important contributions to the program, using ideas from across the mathematical landscape.  In late $2008$, I wrote a long and detailed survey of the state of the art in the program \cite{F11}.  At the time, I was not terribly optimistic about the prognosis for the area.  The most interesting questions seemed
inaccessible and activity in the area seemed to be slowing.  Much to my surprise,
the last decade has proven to be remarkably fruitful for the program and led to a series of breakthroughs
that make the program now more vibrant than ever.  This includes the progress on Zimmer's conjecture
by Brown, Hurtado and me that Bob so graciously highlights in his own contribution to this volume but also several
other related works that made that development possible.  I will very briefly describe the developments
here and point to the subsections in which they are described.  In writing this article, a major
goal is to point to places where similar arguments are used in different contexts, in order to try to isolate
techniques and ideas that are likely to play a key role in further progress on Zimmer's program and in related areas.
Due to constraints of space and time, detailed arguments will not be included and so this paper may serve more as
``reader's guide" to the area than as an introduction or survey.

A recurring theme in this article will be surprising connections and developments.  I start the
narrative with some developments concerning rigidity of Anosov $\Z^d$ actions that end in a proof
of a conjecture of Katok and Spatzier by Rodriguez Hertz and Wang, using some prior results by
Kalinin, Spatzier and myself.  While there has long been some interaction between the Zimmer program and the rigidity of hyperbolic actions of abelian groups, in the past the main successes had concerned results on local rigidity of actions that were a priori hyperbolic or global rigidity results that required very strong assumptions.  Rapid
progress on rigidity of abelian group actions in my work with Kalinin and Spatzier and the work of
Rodriguez Hertz and Wang described in Section \ref{sec:anosovZd} changed the landscape and provided
strong enough tools for $\Z^d$ actions that it made sense to see what would happen for lattice actions.
Spatzier and I encouraged Rodriguez Hertz and Wang to do this, in my case mistakenly believing this
would be more or less a corollary of their theorems.  This led to the breakthrough paper of Brown, Rodriguez Hertz
and Wang on higher rank Anosov actions in which the key philosophy of {\em non-resonance implies invariance} was introduced, see Section \ref{sec:anosovlattice}.  The second application of this philosophy led to a better
understanding of invariant measures for low dimensional actions of lattices in another paper by Brown, Rodriguez Hertz and Wang.  I think it is safe to say that those authors did not appreciate the full importance of what they had done.   Key ingredients throughout these developments are derive from much closer connections
to measure rigidity.  This occurs both in the proofs of the result of Brown, Rodriguez Hertz and Wang and in terms of the
use of dynamics of unipotent flows and particularly Ratner's theorems in the work on Zimmer's conjecture.  All of these developments are described in Section \ref{sec:zimmer}.

In the middle of that narrative an important additional piece of the puzzle was provided by Hurtado's
remarkable paper on the Burnside problem, which formalized a key notion of {\em subexponential
growth of derivatives}.  This notion is really the natural notion of uniform non-hyperbolicity of a group action
and is certainly implicit earlier, but making it explicit and exhibiting as a fulcrum or turning point in proofs
was a key conceptual development, see Section \ref{sec:hurtado}.

Lastly another important piece of the puzzle predates all of these developments: the strong property (T) of
Vincent Lafforgue.  While this notion vaguely resembles one contained in a paper I wrote with Margulis and also some
ideas in a paper of Gromov, the precise notion and its proof are profound and original developments. This is described in Section \ref{sec:strongT}.  Lafforgue developed this notion primarily as an obstruction to
proving the Baum-Connes conjecture and also as a route to constructing superexpanders, but just as property (T) before
it, strong property (T) turns out to be robustly applicable.

This paper is a short and very personal history of these mathematical developments. These developments all come in the context of the Zimmer program and it is a pleasure to describe them here in a volume selecting some of Zimmer's papers.  Following Bob's example in his new essay for this volume, the style here will be one of personal narrative in a manner that is somewhat unusual in modern mathematical exposition.  I want to both explain some mathematical ideas and explain the history of how they arose, combined and lead to exciting new developments.  All of these developments owe a tremendous amount to Bob's insights, theorems and questions.  They also involve the work of several people, each of whom would surely tell this story differently.  I'd like to take the opportunity to thank them all for their roles in these developments including the pleasure I've had in writing papers with most of them.  They are, in alphabetical order, Aaron Brown, Sebastian Hurtado, Boris Kalinin, Federico Rodriguez Hertz, Ralf Spatzier and Zhiren Wang.

\subsection*{Acknowledgments} Thanks to Marc Burger, Manfred Einseidler, Alex Eskin, Simion Filip, Katie Mann, Federico Rodriguez Hertz, and Bob Zimmer for useful comments and corrections.  Particular thanks to Aaron Brown and Ralf Spatzier for thorough readings and copious remarks and corrections.

\section{Pre-history}
\label{sec:prehistory}

I'm going to start the story where it started for me personally.  This was at an AIM workshop organized
by Lindenstrauss, Katok and Spatzier on ``Emerging directions in measure rigidity" in 2004.  From my point
point of view, the start of this long conversation with many people is a remark Ralf Spatzier made during that
problem session.  Either Katok or Spatzier was discussing their conjecture that every {\em genuinely higher rank}
Anosov $\Z^d$ action was smoothly conjugate to an affine Anosov action on an infranilmanifold.  As was usual at the time, discussion turned to the related conjecture, usually attributed to Smale or Franks, that an Anosov diffeomorphisms (or $\Z$ actions) is topologically conjugate to an affine diffeomorphism on a nilmanifold. The conjecture of Katok and Spatzier is motivated by their work on local rigidity and cocycle rigidity for higher rank Anosov actions \cite{MR1307298, MR1632177}. It is well known in the case $d=1$ that one cannot expect smooth conjugacy since perturbations of affine actions can be seen to no longer have constant derivative on the periodic orbits and so are not smoothly conjugate to the original action.  Spatzier reminded the audience of a deeper obstruction from work of Farrell and Jones, the existence of Anosov diffeomorphisms on ``exotic tori", manifolds homeomorphic but not diffeomorphic to tori \cite{MR508890}.  Farrell and Gogolyev give a more modern and general construction in \cite{MR2928077}.  Exotic tori are constructed by taking the connected sum of an exotic sphere and a torus.  The highly non-trivial part of the papers just mentioned is showing that one can take an Anosov diffeomorphism and first restrict it to the complement of a neighborhood of a fixed point, connect sum in an exotic sphere and then extend the Anosov diffeomorphism back over the whole manifold.  Some time later, I realized that no existing result prevented the construction of similar examples of $\Z^d$ Anosov actions
on exotic tori.  I even asked both Katok and Spatzier about this and their responses led me to believe that no one had thought seriously about the possibility.  I then spent some time trying to construct examples, talking about the problem with Chris Connell, Tom Farrell and Shmuel Weinberger at various moments.  The obstruction is a simple one:
if you follow the Farrell-Jones construction, it is completely unclear if the resulting diffeomorphism has non-trivial
centralizer.  Given an exotic Anosov diffeomorphism $f$ on a torus, there are diffeomorphisms $g$ commuting with the $f$ at the level of the homotopy data, so as automorphisms of $\pi_1(\T^n)=\Z^n$.  But it
could easily happen that the commutator of the diffeomorphisms $f$ and $g$, while homotopic to the identity, is non-trivial and the group generated by $f$ and $g$ in $\Diff(\T^n)$ might still be free or at least non-abelian.  Repeated conversation with Spatzier about this issue led to our beginning to work together and with Kalinin and in the not too long run to the resolution of the Katok-Spatzier conjecture for actions on infranilmanifolds by Rodriguez Hertz and Wang, all described in the next section of this article.

\section{Rigidity of Anosov $\Z^d$ actions}
\label{sec:anosovZd}

At the time of the conversations above, the best result concerning Anosov $\Z^d$ actions on tori was
by Federico Rodriguez Hertz.  This paper first appeared as a preprint in 2001 but was only revised and
published in 2007 \cite{MR2318497}.  I must admit I didn't understand this paper well for several years
after that, but it was clear at the time that the dynamical conditions needed were quite restrictive
including needing low dimensional dynamical foliations and the lack of actions on exotic tori in that
context was unsurprising. In the interim between that paper and my work with Kalinin and Spatzier, most of the work on the conjecture of Katok and Spatzier focused on the much harder and still open case of actions on general manifolds. The best work on this topic is by Kalinin and Sadovskya, building on earlier work of Kalinin and Spatzier \cite{MR2322492, MR2372620, MR2240907}. Before continuing, I recall the formal definition of Anosov diffeomorphism
and Anosov action.

Let $a$ be a diffeomorphism of a compact manifold $M$.
We recall that  $a$ is  {\em Anosov} if there exist a continuous $a$-invariant
decomposition  of the tangent bundle $TM=E^s_a \oplus E^u_a$ and constants
$K>0$, $\lambda>0$ such that for all $n\in \mathbb N$
\begin{equation} \label{anosov}
  \begin{aligned}
 \| Da^n(v) \| \,\leq\, K e^{- \lambda n} \| v \|&
     \;\text{ for all }\,v \in E^s_a, \\
 \| Da^{-n}(v) \| \,\leq\, K e^{- \lambda n}\| v \|&
     \;\text{ for all }\,v \in E^u_a.
 \end{aligned}
 \end{equation}
The distributions $E_a^s$ and $E_a^u$ are called the {\em stable} and {\em unstable}
distributions of  $a$. Given an action of a group $\Gamma$ on a compact manifold $M$
via diffeomorphisms we call the action {\em Anosov} if there is an element $\gamma$ in $\Gamma$
such that $\alpha(\gamma)$ is an Anosov diffeomorphism.

Progress on Anosov $\Z^k$ actions on tori and nilmanifolds began again in two papers I wrote with Boris Kalinin and Ralf Spatzier. The starting point for essentially all work on this problem in the context of action on infranilmanifolds
are results of Franks and Manning which show that any Anosov diffeomorphism is topologically conjugate to an affine Anosov diffeomorphism \cite{MR0253352, MR0343317}.  Given an Anosov diffeomorphism $f$ of a torus $\T^n$, we write $f_*$ for the induced linear action of $f$ on $\pi_1(\T^n)=\Z^n$ and note that this defines a linear toral automorphism $f_0$. Franks showed that for any Anosov diffeomorphism on a torus, there is a homemomorphism $\phi$, conjugating an Anosov diffeormorphism $f$ to the linearization $f_0$ \cite{MR0253352}.  This was generalized by Manning to the case of nilmanifolds \cite{MR0343317}.  Since the conjugacy $\phi$ is essentially unique in its homotopy class, it follows that the conjugacy is also a conjugacy for any diffeomorphism commuting with $f$ and so that any Anosov $\Z^d$ action on an infranilmanifold is topologically conjugate to an affine action.  The remaining problem is simply to improve the regularity of the conjugacy.  For this exposition, we restrict attention to the case of actions on tori, where ideas are easier to explain.    The affine model action is always linear on a finite index subgroup.  Call a linear $\Z^d$-action  {\em irreducible} if it does not split rationally  into non-trivial factors.  We can always split a linear action into a product of irreducible factors.
We call a general Anosov action by $\Z^d$
on a $\T^n$ (or more generally a nilmanifold) {\em irreducible} if it is topologically conjugate to an affine action which is irreducible.  A $\Z^k$-action $\alpha$ on $\T^n$  induces a linear action $\alpha _0$ on homology called the {\em linearization} of $\alpha$ and the affine conjugate of an Anosov action is equal to $\alpha_0$ on a subgroup of finite index. The  logarithms of the moduli of the eigenvalues define linear maps $\lambda _i: \Z^k \mapsto \R$ which extend to $\R^k$.  A {\em Weyl chamber } of $\alpha _0$ is a connected component of $\R^k - \cup _i \ker \lambda _i$.  The key  hypothesis in the work with Kalinin and Spatzier is the presence of an Anosov diffeomorphism in each Weyl  chamber and a special case of our results in \cite{MR2776843, MR2983009} is

 \begin{theorem}\label{theorem:katokspatzier}  Suppose $\alpha $ is an irreducible $\ci$-action of $\Z ^k$, $k \geq 2$ on a torus $\mathbb T^n$.  Further assume there is an Anosov element for $\alpha$ in each Weyl chamber of $\alpha _0$. Then $\alpha$ is $\ci$-conjugate to an affine action with linear part $\alpha_0$.
\end{theorem}

\noindent The key point of the assumption of having an Anosov in each Weyl chamber is the following.  Just using Oseledec theorem and Pesin theory, one has that the action $\alpha$ is non-uniformly partially hyperbolic in a particularly nice way.  The presence of enough Anosov diffeomorphisms forces all invariant measure to have proportional Lyaponuv exponents and therefore forces essentially all elements of the action to be uniformly hyperbolic in the strongest possible sense. The fact that controlling all invariant measures turns non-uniform estimates obtained from Oseledec theorem into uniform estimates is an old one in dynamical systems, is also used in Rodriguez-Hertz's earlier paper and is also used in a surprising way in the proof of Zimmer's conjecture.  In particular, in this context, it shows that if there are elements of $\Z^d$ in Weyl chamber walls, then restricted to the corresponding {\em central foliation} one has {\em subexponential growth of derivatives}.  This viewpoint is a bit ahistorical as the notion of subexponential growth of derivatives for a group action was only formalized later by Hurtado in \cite{Hurtado}.  We will discuss the idea in detail in the less technical setting of that article in Section \ref{sec:hurtado} below.

In the first of our papers, this observation is combined with ideas concerning normal forms for group actions
and an examination of holonomies to produce a result for so-called {\em totally non-symplectic actions} \cite{MR2776843}.  While the techniques used in this paper have not had any further applications that I know of, it was around this time that I started giving talks emphasizing the close connections to homogeneous dynamics that had
last been clear in early work of Katok and Spatzier  \cite{MR1632177, MR1406432}, since our holonomy arguments were inspired by the use of unipotents in work of both Ratner and Lindenstrauss on rigidity of invariant measures.  This connection to measure rigidity turns out to have been deep and fruitful as we will see below.

In our second paper \cite{MR2983009}, we took an entirely different approach that turns out to also have some remarkable similarities to the work on Zimmer's conjecture.  Namely we wrote the conjugacy, at least projected
onto certain dynamically defined submanifolds for the action, as a solution to a cohomological equation.  We then combined the subexponential growth of derivatives with exponential decay of matrix coefficients
to obtain greater regularity of the solution to the cohomological equation.  This argument is quite complicated
as one needs to work along various foliations and so the regularity one obtains is best described
in terms of wavefront sets, which then allow one to patch the regularity together globally with
arguments that are fairly standard in PDE.  The remarkable thing is that this argument is quite
similar to the last step in the proof of Zimmer's conjecture, where we use strong property (T) in conjunction
with subexponential growth of derivatives (in all directions) to find an invariant metric.  The proof of
strong property (T) depends essentially on exponential decay of matrix coefficients, so the similarity
in the arguments appears to be quite deep.  There is, of course, an important point of contrast. In our
setting the solution to the cohomological equation is given.  Lafforgue proves additional estimates of decay of
certain twisted matrix coefficients in order to prove the solution exists.

The paper of Rodriguez Hertz and Wang completes the proof of the Katok-Spatzier conjecture
for higher rank abelian Anosov actions on nilmanifolds and tori \cite{MR3260859}.  They proceed by showing that the
presence of a single Anosov element in a $\Z^d$ action on a nilmanifold or torus implies the existence of an Anosov element in every Weyl chamber and therefore verify the hypothesis of Theorem \ref{theorem:katokspatzier} above. It is relatively easy to see that the presence of the single Anosov element forces the entire $\Z^d$ action to be non-uniformly hyperbolic in a strong sense and so the key contribution is showing that this non-uniform hyperbolicity is uniform.  This is a theme that recurs in both their work with Brown on rigidity of Anosov actions of higher rank lattices on tori and nilmanifolds and in the proof of Zimmer's conjecture. I expect this them will continue to recur in further progress in the program.

 The first step in their argument is to show that the conjugacy is smooth along certain dynamically defined
foliations, namely the foliations along which the map is switching from expansion to contraction as one changes Weyl chambers.  The proof of  this fact follows quite closely arguments in \cite{MR2983009}. The main part of the paper is then to show that this smooth map along dynamical foliations, in addition to being smooth, is non-singular.  That this suffices is shown using a theorem of Man\~{e}, which shows that any quasi-Anosov diffeomorphism is in fact Anosov.  Man\~{e}'s work seems a fertile source of ideas for understanding further connections between non-uniform hyperbolicity and uniform hyperbolicity.

The basic idea for the hardest part of the paper is that along the singular set, the conjugacy collapses volume.  Rodriguez Hertz and Wang derive a contradiction to this fact by using a measure rigidity type argument to show that this collapse of volume is impossible.  This is by far the most delicate part of the paper and uses that strong stable manifolds are Lipschitz inside of stable manifolds.  The measure rigidity result is particularly difficult since the measure considered is an invariant measure supported on the (closed) set of singular points. This application of measure rigidity ideas in the non-linear settings and in particular using them to improve non-uniform hyperbolicity to uniform hyperbolicity marks another instance of an important new trend in Zimmer's program.  We will continue to see more instances of this below.

The reader far from hyperbolic dynamics may wonder why the focus on converting non-uniform hyperbolicity to uniform hyperbolicity when the real breakthrough comes on Zimmer's conjecture, which concerns actions with no hyperbolicity.  One motivation for the conjecture was Zimmer's cocycle super-rigidity theorem, which produced a measurable invariant metric and showed all Lyapunov exponents were zero.  While the program traditionally was framed as an attempt to improve regularity of this metric, the breakthrough comes from a change of perspective.  Rather than improving regularity of the metric directly, one improves estimates on dynamics of the derivative from non-uniform to uniform.

In closing, it is an old conjecture often attributed to either Smale or Franks that all Anosov diffeomorphisms are conjugate to affine maps on nilmanifolds.  In the final comment of \cite{MR1754775}, Margulis draws an intriguing parallel between now resolved cases of Zimmer's conjecture, open problems in the Zimmer program and this
conjecture on Anosov diffeomorphisms.  Resolving this conjecture, combined with the work just described would provide a smooth classification of Anosov actions of higher rank abelian groups.

\section{Rigidity of Anosov lattice actions}
\label{sec:anosovlattice}

The next major breakthrough concerned actions of higher rank lattices, not just higher rank abelian groups,  and so was  truly part of Zimmer's program. It is still deeply embedded in hyperbolic dynamics as it concerns Anosov actions.  The paper of Brown, Rodriguez Hertz and Wang contains not one but two major results each containing ideas which are relevant to future directions \cite{BRHW1}.  We recall that for any group $\Gamma$ an action is {\em Anosov} if some individual element $\gamma$ in the group acts as an Anosov diffeomorphism.

In the first part of the paper, they give a proof that any Anosov action of a higher rank lattice $\Gamma$ on a torus or nilmanifold $M$ is continuously conjugate to an an affine action subject to the condition that the action lifts to the universal cover.  For the case of actions preserving a measure of full support, this reproves a result of Margulis and Qian \cite{MQ01}, but the proof by Brown, Rodriguez Hertz and Wang is different and does not use Zimmer's cocycle super-rigidity theorem.  The philosophy behind this argument is an important new development that is also applied in the proof of Zimmer's conjecture.  This philosophy develops a link between the finite dimensional representation theory of the group $G$ and the possible actions of $G$ and $\Gamma$.  A link between the two appeared already in Zimmer's cocycle superrigidity theorem. Zimmer showed that, for smooth actions of $G$ or $\Gamma$, the derivative of measure preserving action is always essentially described, at least measurably, by a finite dimensional representation of $G$. The presence of an invariant measure is a necessary hypothesis for Zimmer's result.  The philosophy of ``non-resonance implies invariance" allows one to extend the connection between group actions and finite dimensional representations to the setting where there are, a priori, no invariant measures.

To employ this philosophy to actions of  a lattice $\Gamma$ one always need to pass to the induced $G$ action on $(G \times M)/\Gamma$. This allows one to use the structure of $G$, namely the root data associated to a choice of Cartan subalgebra. To explain this philosophy better, I recall some basic facts.  The Cartan subgroup $A$ of $G$ is the largest subgroup diagonalizable over $\R$, the Cartan subalgebra $\mathfrak a$ is it's Lie algebra.  It is known since the work of \'{E}lie Cartan that a finite dimensional linear representation $\rho$ of $G$ is completely determined by linear functionals on $\mathfrak a$ that arise as generalized eigenvalues of the restriction of $\rho$ to $A$.  Here we use that there is always a simultaneous eigenspace decomposition for groups of commuting symmetric matrices and that this makes the eigenvalues into linear functionals.   These linear functionals are referred to as the {\em weights} of the representation.  For the adjoint representation of $G$ on it's own Lie algebra, the weights are given the special name of {\em roots}.  Corresponding to each root $\beta$ there is a unipotent subgroup $G_{\beta} <G$ called a {\em root subgroup} and it is well known that ``large enough" collections of root subgroups generate $G$.  Two linear functionals are called {\em resonant} if one is a positive multiple of the other.  Abstractly, given a $G$-action and an $A$-invariant  object $O$, one may try to associate to $O$ a class of linear functionals $\Omega$. {\em Non-resonance implies invariance} is the observation  that, given any root $\beta$ of $G$ that is not resonant to an element of $\Omega$, the object $O$ will automatically be invariant under the unipotent root group $G^{\beta}$.  If one can find enough such non-resonant roots, the object $O$ is automatically $G$-invariant.

The philosophy  that``non-resonance implies invariance" is perhaps more transparent in it's application to invariant measures in Section \ref{sec:zimmer}, but I will describe it here first. To describe this philosophy in action here, we need to develop the picture a bit further.   As discussed in the last section, Franks and Manning produced a conjugacy $\phi$ conjugating an Anosov diffeormorphism $f$ to the linearization $f_0$. Note that the action on fundamental group gives a linear representation of $\Gamma$ into $Aut(\pi_1(M))$ restricting to $f_0$ on the Anosov diffeomorphism $f$.  Using this $\Gamma$ action on homology, the structure of $M$ as a nilmanifold and Margulis' superrigidity theorem, one defines an affine $G$ action on $(G \times M)/\Gamma$. As discussed above, the map $\phi$ is sufficiently unique that it also conjugates the centralizer of the Anosov diffeomorphism to linear maps.  Letting $A$ be the maximal Cartan subgroup of $G$, a variant of that argument produces a conjugacy $\phi$ between the restriction to $A$ of  the original $G$-action  and the affine action of $A$  on $(G \times M)/\Gamma$. This affine action is defined in terms of a linear representation of $A$ and therefore gives rise to a collection  $\Omega$ of linear functionals on $A$. More concretely, since $A$ acts affinely and $(G \times M)/\Gamma$ is parallelizable, the derivative is constant and these linear functionals are exactly generalized eigenvalues of the derivative.   So this semiconjugacy to a linear action gives rise to a new collection $\Omega$ of linear functionals on $A$ that one can compare to the roots of $G$.  Brown, Rodriguez Hertz and Wang show that for any root $\beta$ that is not resonant to any linear functional $\alpha$ in $\Omega$, the map $\phi$ is $G_{\beta}$ invariant. One key point for the proof is that the linear functionals in $\Omega$ associated to Anosov actions cannot be resonant to roots for purely algebraic reasons. The other key point is the existence of elements of $A$ that commute with $G_{\beta}$. This shows that $\phi$ is $G$ equivariant, as desired.

The second half of the paper proves that this continuous conjugacy is smooth.  The argument surprisingly goes back to an old argument in a paper of Katok, Lewis and Zimmer on actions of $\Sl(n,\Z)$ on $\T^n$ \cite{KLZ}.  Both there and here the goal is to show that the Lyapunov exponents coming out of cocycle superrigidity can be improved to give uniform estimates on derivatives.  In this paper, the key point is to control the Lyapunov exponents uniformly along all invariant measures, using that Ratner's theorem classifies invariant measures on the linear side and so also on the non-linear side by conjugacy.  This combination of cocycle superrigidity with control on exponents in all invariant
measures also foreshadows elements of the proof of Zimmer's conjecture and will certainly appear in future developments as well.

\section{Hurtado's work on the Burnside problem for diffeomorphisms}
\label{sec:hurtado}

A key impetus for Brown, Hurtado and I to work together on Zimmer's conjecture was Hurtado's paper proving a Burnside type theorem for subgroups of $\Diff(S^2,\omega)$, i.e for volume preserving diffeomorphisms of the sphere in dimension $2$ \cite{Hurtado}.  Hurtado proves that any finitely generated subgroup $\Gamma$ of $\Diff(S^2, \omega)$ that consists entirely of torsion elements is finite. Since the conjecture was known for all other surfaces, the focus on $S^2$ is natural. The requirement of an invariant volume form is an artifact of the proof.   Farb and Ghys have each independently conjectured finiteness of all finitely generated torsion subgroups of the diffeomorphism group of any compact manifold.

Hurtado formalizes a key notion of subexponential growth of derivatives.   Let $\Gamma$ be a finitely generated group.  Let $\wl\colon \Gamma\to \N$ denote the word-length function relative to some fixed finite set of generators $F$. Let $\alpha \colon \Gamma \to \diff^1(M)$ be an  action of $\Gamma$  on a compact manifold $M$ by $C^1$ diffeomorphisms.   We say the action $\alpha$ has \emph{uniform subexponential growth of derivatives} if for all $\epsilon>0$ there is a $C_{\epsilon}$ such that for all $\gamma \in \Gamma$ we have
$$\|D\alpha(\gamma)\|\le C_{\epsilon}e^{\epsilon \wl(\gamma)}$$
where $\|D\alpha(\gamma)\|= \sup_{x \in M} \|D_x\alpha(\gamma)\|$.

In Zimmer's own work, several notions of slow growth of derivatives arose and showing slow enough growth
of derivatives was clearly the key to the conjecture from many points of view. I will discuss this more at
the end of Section \ref{sec:zimmer} below.  In the context of hyperbolic dynamics the notion of subexponential growth here is the correct uniform analogue of an action having all zero Lyapunov exponents. Even though this idea has been used before (even implicitly in my work with Kalinin and Spatzier), Hurtado’s paper seems to be the first place it is formalized for group actions and its formulation proved extremely useful for our subsequent work on Zimmer’s conjecture.  For a single diffeomorphism it is an easy classical fact that having zero Lyapunov exponents for all invariant measures implies subexponential growth of derivatives.  I state this result formally here as it, and Hurtado's adaptation of it, were starting points for the work on Zimmer's conjecture.

\begin{proposition}
\label{prop:easydynamics}
Let $M$ be a compact manifold and let $\Z$ act smoothly on $M$.  Then the $\Z$ action has subexponential growth
of derivatives if and only if every $\Z$ invariant measure on $M$ has zero Lyaponuv exponents.
\end{proposition}

\noindent To prove the non-trivial implication in the proposition, one takes sequences of orbit segments witnessing the failure of subexponential growth of derivatives, views these as defining measures and shows the weak$*$ limits of these measures have positive Lyapunov exponents.

Hurtado uses an analogue of this result for the $\Z$ action built by taking the shift action of $\Z$ on the space  $\Omega=F^{\Z}$ and looking at the standard skew product action of $\Z$ on $\Omega \times S^2$.  An extension of the classical argument says that either the $\Gamma$ action has subexponential growth of derivatives or there is a $\Z$ invariant measure on $\Omega \times S^2$ with a positive exponent
for the derivative along $S^2$.  Since Hurtado assumes an invariant volume form, the presence of one non-zero exponent
forces both exponents along $S^2$ to be non-zero.  At this point, Hurtado invokes a classical result of Katok on hyperbolic dynamics of diffeomorphisms with no zero exponents \cite{MR573822}.  To use this theorem, he needs to embed  $\Omega \times S^2$ in a manifold and construct a hyperbolic diffeomorphism $f$ of the manifold extending the shift action on $\Omega \times S^2$.  Katok's theorem then provides a hyperbolic fixed point for $f$ which in this construction gives a hyperbolic periodic point for some element of $\Gamma$ on $S^2$.  Since all elements of $\Gamma$ are finite order, this is a contradiction.

Hurtado's paper contains many other intriguing ideas once this result is in place. In fact, he uses the idea of averaging a metric over an action with subexponential growth of derivatives to obtain an invariant metric. Here this only works directly under the additional assumption that the group is amenable, but this too foreshadows elements of the proof of Zimmer's conjecture. Experts appear to have known that subexponential growth of derivatives should be relevant to the Burnside problem before Hurtado's work. The approach is at least partially inspired by Kalinin's matrix valued Liv\v{s}ic theorem \cite{MR2776369} and the diffeomorphism valued Liv\v{s}ic theorem of Kocsard and Potrie \cite{MR3471936}, but Hurtado's success in implementing this style of proof was a major impetus for our work on Zimmer's conjecture.  The connection to Liv\v{s}ic type theorems is also related to another paper of Hurtado's  where he produces a surprising example \cite{HurtadoII}.  He produces an action of a free group
where every element is conjugate to an isometry but the full action has exponential growth of derivatives.

\section{Strong property (T)}
\label{sec:strongT}

Another key ingredient in the proof of Zimmer's conjecture that was developed much earlier is
Lafforgue's notion of strong property $(T)$ introduced in 2007 by Lafforgue as an obstruction to certain strategies for proving the Baum-Connes conjecture \cite{MR2423763}.  The original definition is made in terms of existence of certain kinds of projections in a certain completion of the group algebra. Instead, we recall a variant that does not involve operator algebras.
Given a group $\Gamma$ and a finite or compact generating set $S$, we let $\wl$ be the word length on $\Gamma$.
In what follows $X$ will denote a Banach space and $B(X)$ will denote the bounded operators on $X$. We always consider the operator norm topology on $B(X)$ and we always mean the operator norm when we write $\|T\|$ for $T \in B(X)$.

\begin{definition}
\label{defn:subexpnorm}
Let $\Gamma$ be a group with a length function $\wl$, $X$ a Banach space and $\pi\colon \Gamma \rightarrow B(X)$.
Given $\epsilon>0$, we say $\pi$ has \emph{$\epsilon$-subexponential norm growth} if there exists a constant $L$ such that
$\|\pi(\gamma)\| \leq Le^{\epsilon\wl(\gamma)}$ for all $\gamma \in \Gamma$.
We say $\pi$ has subexponential norm growth if it has $\epsilon$-subexponential norm growth for all
$\epsilon>0$.
\end{definition}

Here we will focus only on the case of actions on Hilbert spaces, though the robustness of strong property $(T)$ outside of that class is important for many applications and is discussed briefly below.

\begin{definition}
\label{defn:strongT}
A group $\Gamma$ has \emph{strong property (T)} if there exists a constant $t>0$ and a sequence  of measures $\mu_n$ supported in the balls $B(n) = \{\gamma \in \Gamma \mid  l(\gamma) \leq n\}$   in $\Gamma$ such that: for any representation $\pi\colon  \Gamma \rightarrow B(X)$ with $t$-subexponential norm growth and $X$ a Hilbert space,
the operators $\pi(\mu_n)$ converge exponentially quickly to a projection onto the space of invariant vectors.  That is,  there exists  $0<\lambda<1$ (independent of $\pi$),  a projection $P \in B(X)$ onto the space of $\Gamma$-invariant vectors, and an $n_0\in \N$ such that $\|\pi(\mu_n)-P\| <  \lambda^n$ for all $n\ge n_0$.
\end{definition}

The original definition of Lafforgue is equivalent not to this but to a similar statement with signed measures
in place of measures \cite{delasalle}.  All known proofs that a group has strong property $(T)$ produce positive measures and this definition seems more useful, if harder to state in the purely operator theoretic language \cite{delasalle}.  If one makes similar
definitions but assumes instead that $L=1$ in Definition \ref{defn:subexpnorm} then the resulting notion is equivalent to property $(T)$ by my work with Margulis in \cite{MR2198325}.  One might say that in that work we insist on {\em immediate} subexponential norm growth while Lafforgue allows the weaker condition of {\em eventual} subexponential norm growth.  The difference is quite profound. Lafforgue prove that no hyperbolic groups can have strong property $(T)$ while many (even most in certain random senses) of them have property $(T)$.  The difference is also profound for applications. One should only expect immediate subexponential growth of derivatives when perturbing isometric actions  as in \cite{MR2198325}.

By work of Lafforgue, de la Salle and de la Salle-de Laat, strong property $(T)$ is known for the groups that interest us here: \cite{MR2423763, MR3407190, delaSallenonuniform}

\begin{theorem}
\label{theorem:strongT}
Let $G$ be a semisimple Lie group all of whose simple factors have real rank at least two and let $\Gamma<G$
be a lattice.  Then $G$ and $\Gamma$ have strong property $(T)$.
\end{theorem}

\noindent Combined with additional work of Liao, these results also establish strong property $(T)$ for lattices
in higher rank simple groups over other local fields \cite{MR3190138}.

While the notion just described suffices to prove Zimmer's conjecture for smooth actions,  to obtain Zimmer's conjecture in lower regularity, one needs to consider non-Hilbertian function spaces.  In this context, it suffices to consider the $\theta$-Hilbertian Banach spaces.  These are subspaces of spaces obtained by complex interpolation with Hilbert spaces.  Most reasonable function spaces arising in dynamics and geometry are $\theta$-Hilbertian except those which are defined in terms of supremum norms and so have no good convexity properties.  Strong property $(T)$ for higher rank lattices is known to hold for all $\theta$-Hilbertian Banach spaces and even for significantly broader classes of Banach spaces. It is conjectured to hold for any uniformly convex Banach space.

\section{Recent work on Zimmer's conjecture}
\label{sec:zimmer}

In this section, I will focus on the recent breakthrough on Zimmer's conjecture, discussing only the case of actions
of lattices in $\Sl(n,\R)$ for $n>2$ since that simplifies terminology and notation. For an excellent account both of the developments discussed here and of the numerology associated to actions of other groups, the reader should consult Cantat's recent Bourbaki Seminaire article \cite{Cantat}.  Another account of the recent developments on Zimmer's conjecture that is probably more introductory and accessible than the one here is contained in a survey by Brown \cite{BrownSurvey}.

The discussion here will focus on three articles.  First the article of Brown, Rodriguez Hertz and Wang that produces invariant measures for low dimensional actions \cite{AWBFRHZW-latticemeasure}, then the article by Brown, Hurtado and myself which proves Zimmer's conjecture for actions of cocompact lattice \cite{BFH}, at least for $\R$ split classical groups, and lastly the very recent preprint in which Brown, Hurtado and I establish the conjecture for $\Sl(n,\Z)$ for $n>2$ \cite{BFH2}. Throughout this section we fix the notation that $G=\Sl(n, \R)$ and $n>2$ and
$\Gamma <G$ a lattice.

The first results we refer to are those of Brown, Rodriguez Hertz and Wang on rigidity of invariant measures.  A simple version of their main result is

\begin{theorem}
\label{theorem:invariantmeasure}
Assume $\Gamma$ acts smoothly on a compact manifold with dimension less than $n-1$.  Then there is a $\Gamma$ invariant measure on $M$.
\end{theorem}

\noindent Their results are much stronger than this and also find an invariant measure in certain higher dimensions unless the action has a {\em measurable projective factor}, i.e. a measurable quotient given by the $\Gamma$ action $G/Q$ where $Q$ is a parabolic. They also have criteria for which parabolic subgroups can arise in terms of the dimension of $M$.  The theorem has some resemblance to work of Nevo and Zimmer that also produces invariant measures or  measurable projective factors \cite{MR1933077}, but the important difference is a computable criterion for triviality of the measurable projective factor in terms of Lyapunov exponents associated to a measure invariant by the maximal Cartan subgroup. While the techniques in Nevo and Zimmer are a difficult generalization of the proofs of Margulis' projective quotient theorem, the proof here is a less direct rethinking of that proof and can be used to give a proof of the projective factor theorem that is organized somewhat differently than Margulis'.

I want to describe the strategy in some detail as it is another example of the philosophy that non-resonance implies invariance. The reader may want to reread the first paragraph or two describing this philosophy in Section \ref{sec:anosovlattice} before proceeding, at least to recall the basic background and motivation.  Recall that two linear functionals are called resonant if one is a positive multiple of the other. As in Section \ref{sec:anosovlattice} one first  passes to the induced action on $X= (G \times M)/\Gamma$. Taking the minimal parabolic $P$ and using that it is amenable, one finds a $P$ invariant measure $\mu$ and the goal is to prove that $\mu$ is $G$ invariant.  The measure $\mu$ is also clearly invariant under the Cartan subgroup $A$ which is contained in $P$ and so one can try to apply the philosophy that non-resonance implies invariance by associating some linear functionals to the pair $(A, \mu)$.  The linear functionals we consider are the Lyapunov exponents for the $A$ action. More precisely we consider the  Lyapunov exponents for the restriction of the derivative of $A$ action to the subbundle $F$ of $T((G \times M)/\Gamma)$ defined by directions tangent to the $M$ fibers in that bundle over $G/\Gamma$.  We refer to this collection of linear functionals as {\em fiberwise Lyapunov exponents}. In this context \cite[Propsition 5.3]{AWBFRHZW-latticemeasure} shows that, given an $A$ invariant measure on $X$ that projects to Haar measure on $G/\Gamma$, if a root $\beta$ of  $G$ is not resonant with any fiberwise Lyapunov exponent then the measure is invariant by the root subgroup $G_{\beta}$.  The rest of the proof is quite simple.   The stabilizer of $\mu$ contains $P$ which implies the projection of $\mu$ to $G/\Gamma$  is Haar measure, so the proposition just described applies.  The stablizer $G_{\mu}$ of $\mu$ in $G$ is a closed subgroup containing $P$. We also know that $G_{\mu}$  contains the group generated by the $G_{beta}$ for all roots $\beta$ not resonant with any fiberwise Lyapunov exponent. We also know that the number of distinct fiberwise Lyapunov exponents is bounded by the dimension of $M$.  Since any closed subgroup of $G$ containing $P$ is parabolic, $G_{\mu}$ is parabolic.  So either $G_{\mu}=G$ or the number of resonant roots needs to be at least the dimension of $G/Q$ for $Q$ a maximal proper parabolic. This is because given any single root $\beta$ with $G_{\beta} \nless Q$ the group generated by $G_{\beta}$ and $Q$ is $G$.  Our assumption on the dimension of $M$ immediately implies there are not enough fiberwise Lyapunov exponents to produce $\dim(G/Q)$ resonant roots, so $\mu$ is $G$ invariant.

The hard part of the proof is contained in \cite[Propsition 5.3]{AWBFRHZW-latticemeasure}.  This follows an outline that is common in measure rigidity, in that the key tool is computing entropy and using deep work of Ledrappier and Young relating entropy to dimensions of invariant measures \cite{MR743818, MR819556}.  To apply these techniques, they also need to redevelop the Ledrappier-Young theory along coarse Lyapunov foliations for the group action and a more refined form of the Abramov-Rohklin theorem for entropy of skew products \cite{AWB-GLY}.  We remark here that for lattices in $SL(n,\mathbb R)$ one can get away with a less difficult application of the work of Ledrappier-Young that was discovered by Hurtado and is explained by Cantat in \cite{Cantat}, see below for more discussion.

In this context, the main result proven in the work with Brown and Hurtado \cite{BFH} is

\begin{theorem}
\label{theorem:invariantmetric}  Assume $\Gamma$ is a cocompact lattice and
assume $\Gamma$ acts smoothly on a compact manifold $M$ with dimension less than $n-1$.  Then the action factors through a finite quotient of $\Gamma$.  Suppose $\Gamma$ acts smoothly on a compact manifold $M$ of dimension $n-1$
preserving a smooth volume form.  Then again the action factors through a finite quotient.
\end{theorem}

The proof of this theorem again begins by inducing the action to a $G$ action on $(G \times M)/\Gamma$.  Since $\Gamma$ is cocompact it is quite easy to verify that subexponential growth of derivatives for the $\Gamma$ action on $M$ is equivalent to subexponential growth of derivatives for $G$ along the vector bundle $F$ tangent to the $M$ fibers in $(G \times M)/\Gamma$.  The derivative along $G$ is computed by the adjoint action of $G$ on its Lie algebra so has exponential growth.  We denote by $F$ the subbundle of $T(G\times)M)/\Gamma)$ that is tangent to the $M$ fibers and refer to the derivative restricted to $F$ as the {\em fiberwise derivative} and Lyapunov exponents for the derivative restricted to $F$ as {\em fiberwise Lyapunov exponents}.  We would like to control growth of fiberwise derivatives
by controlling fiberwise Lyapunov exponents.  For any $G$ invariant measure, it follows from Zimmer's cocycle superrigidity theorem that all fiberwise Lyapunov exponents are zero.  But $G$ and $\Gamma$ are non-amenable so a priori one does not have many $G$ invariant measures with which to work.

The first key observation in our proof is to recall that we have a Cartan decomposition $G=KAK$ where $A$ is the diagonal matrices and so isomorphic to $\R^{n-1}$ and $K=SO(n)$ is compact.  Using compactness of $K$ to average,
it is obvious that there is a $K$ invariant metric on $(G\times M)/\Gamma$ and so subexponential growth of fiberwise derivatives for $G$ is equivalent to subexponential growth of fiberwise derivatives for $A$.  Using essentially the proof of Proposition \ref{prop:easydynamics}, we can then show that if the $G$ action does not have subexponential growth of fiberwise derivatives, there is an $A$ invariant measure $\mu$ with non-zero fiberwise Lyapunov exponents.  The heart of our proof is to promote this measure to one that is $G$ invariant and so obtain a contradiction to Zimmer's cocycle superrigidity theorem.

The desired $G$ invariance is obtained in two steps.  There is a natural projection
$\pi: (G \times M)/\Gamma \rightarrow G/\Gamma$ and the obstruction to applying the proof of Theorem \ref{theorem:invariantmeasure} directly is that we do not know $\pi_*(\mu)$ is Haar measure or that $\mu$ is $P$ invariant.  We begin by carefully averaging $\mu$ over various unipotent subgroups of $G$.  We need to do this in a manner which

\begin{enumerate}
\item preserves the $A$ invariance of $\mu$,
\item preserves the positive fiberwise Lyapunov exponent of $\mu$ and
\item increases the size of the stabilizer of $\pi_*(\mu)$.
\end{enumerate}

\noindent Doing this requires that we use Ratner's theorems on measures invariant under unipotent flows as well as some additional facts from her proof about stabilizers of measures and a result about uniqueness of averages due to Ratner for one parameter unipotent subgroups and to Shah for averages over higher dimensional unipotent subgroups \cite{MR1135878, MR1262705, MR1291701}.

Once we know that $\pi_*(\mu)$ is Haar, we still cannot argue as in the proof of Theorem \ref{theorem:invariantmeasure}, since we do not know the stabilizer of the invariant measure is parabolic.  However, it turns out there are no closed subgroups of $G$ of dimension less than $n-1$, so for manifolds of dimension less than $n-2$ just counting the number of roots of $G$ whose root groups can be resonant with fiberwise Lyapunov exponents completes the proof.  To finish the volume preserving proof, one uses that the only subgroups of codimension $n-1$
are parabolic and that the roots omitted in one of those cannot be the fiberwise Lyaponuv exponents of a volume preserving action.

Once we've established that derivatives of the $\Gamma$ action grow subexponentially, the proof hinges on strong property $(T)$ but is otherwise quite close to arguments in my paper with Margulis on local rigidity of isometric actions \cite{MR2198325}.  Using a variant of the chain rule for higher derivatives, one verifies that subexponential growth of derivatives implies subexponential norm growth for actions of $\Gamma$ on various Sobolev type spaces of symmetric two forms.  Picking a smooth metric $g$ one then looks at the averages
of these over the $\mu_n$ provided by strong property $(T)$ as in definition \ref{defn:strongT}.  These converge
to a $\Gamma$ invariant symmetric two form $g_0$.  Using that the convergence to the invariant form is exponential while the derivatives are subexponential, one checks that $g_0$ is positive definite and so a metric. The Sobolev embedding theorems allow one to see $g_0$ is smooth and the rest of the proof is standard.

I want to make a short remark on a historical quirk. The proof of Theorem \ref{theorem:invariantmetric} in dimension less than $n-1$ can be simplified substantially as explained in \cite[Section 8.3]{Cantat}. Though it is not remarked explicitly there, this simplification also applies to volume preserving actions in dimension $n-1$. This also gives a much simpler proof of Theorem \ref{theorem:invariantmeasure} in dimension less than $n-1$. This simplified version doesn't really use the philosophy of non-resonance implying invariance, and doesn't flow as naturally from the history of ideas discussed in this essay.  It is, in fact, somewhat surprising that this easier version, particularly as a proof of Theorem \ref{theorem:invariantmeasure} was not discovered much earlier.  To me this justifies my choice of writing this essay as a history of ideas that indicates how we reached this point, rather than writing an article that simply explains the mathematics that exists at the end of the trajectory.  Mathematical ideas may be logical but their development is an idiosyncratic human activity and the failure to proceed along the shortest line to a goal is often fruitful in surprising ways.  I should also remark that the easier proofs produce much less sharp results for most simple Lie groups.

I will discuss briefly the issues that arise in the case of non-uniform lattices.  This is discussed in much more detail in the introduction to the recent preprint with Brown and Hurtado proving Zimmer's conjecture for actions of $\Sl(m,\Z)$ where $m>2$ \cite{BFH2}.  In this case $(G\times M)/\Gamma$ is not compact, so the space of probability measures on it is not weak$*$ compact. This presents difficulties both for producing invariant measures with non-zero Lyapunov exponents and then averaging them to improve their projections.  Whether we are taking a given measure and translating it or taking a sequence of empirical measures, it is unclear that the sequences of measures we consider have limits.   At the step of averaging over unipotents, this issue is easily controlled using results on quantatitve non-divergence of unipotent orbits from work of Kleinbock-Margulis \cite{MR1652916, MR2648694}.  It is worth remarking that the first results on non-divergence of unipotent orbits were part of Margulis' dynamical proof of arithmeticity of non-uniform lattices \cite{MR0422499, MR0463353, MR0463354}, a proof that preceded his superrigidity theorem by a few years. By the late 70s there were proofs of arithmeticity for non-uniform lattices via superrigidity that did not pass through homogeneous dynamics, but I find it quite striking that to prove Zimmer's conjecture we needed ideas from both proofs of arithmeticity.

The problem of constructing $A$ invariant measures to begin with is more serious and seems related to issues and conjectures regarding the set of divergent $A$ orbits in $G/\Gamma$ as discussed in the paper of Kadayev, Kleinbock, Margulis and Lindenstrauss \cite{KKLM}.  It is not clear that it is possible to resolve this issue directly, even modulo several important conjectures in homogeneous dynamics, but this is not the only obstruction in this setting. Several other problems arise due to the properties of the {\em return cocycle} $\beta: G \times G/\Gamma \rightarrow \Gamma$ in this context.  The first of these is that subexponential growth of derivatives for the $\Gamma$ action is not equivalent to subexponential growth of derivatives for the induced $G$ action on $(G \times M)/\Gamma$. Fixing a basepoint $x_0$ in $G/\Gamma$, if $x_n$ are points in $G/\Gamma$ with $d(x_n,x_0)=n$, then for any fixed $g$, the size of $\beta(g,x_n)$ will grow linearly in $n$.  Because of this, subexponential growth of fiber derivatives for the induced action is equivalent to having a $\Gamma$ invariant metric on $M$.  Despite this, we still manage to construct a proof using the induced action.  Another key difficulty for all approaches is that we are not able to control the ``images" of the cocycle $\beta\colon G \times G/\Gamma \rightarrow \Gamma$ even for $\Sl(n,\Z)$. To understand this remark better, consider first the case where $G=\Sl(2,\R)$ and $\Gamma =\Sl(2,\Z)$. If we take a one-parameter subgroup $c(t)<\Sl(2,\R)$ and take the trajectory $c(t)x$ for $t$ in some interval $[0,T]$ and assume the entire trajectory on $G/\Gamma$ lies deep enough in the cusp, then $\beta(a(t),x)$ is necessarily unipotent for all $t$ in $[0,T]$.  No similar statement is true for $G=\Sl(m,\R)$ and $\Gamma=\Sl(m,\Z)$.  In fact analogous statements are true if and only if $\Gamma$ has $\Q$-rank one. This is closely related to the fact that higher $\Q$-rank locally symmetric spaces are $1$-connected at infinity.  This forces us to ``factor" the action into actions of rank-one subgroups in order to  control the growth of derivatives.  In the case of $\Sl(m,\Z)$ we prove subexponential growth of derivatives for the $m^2$ canonical copies of $\Sl(2,\Z)$ in $\Sl(m,\Z)$ given by choices of pairs of coordinates and then use the result of Lubotzky, Mozes and Raghunathan that  $\Sl(m,\Z)$ is boundedly generated by these to finish the proof \cite{MR1244421}.  We are currently working on a proof for general non-uniform lattices jointly with Dave Witte Morris.  One step is to find an analogue of this bounded generation statement in general. While this is how Lubotzky, Mozes and Raghunathan prove their results on non distortion of $\Sl(m, \Z)$ in $\Sl(m,\R)$ it is not how their proof proceeds for general higher rank lattices \cite{MR1828742}.  Our proof in the general case also makes use of some more recent results in homogeneous dynamics.

Before closing this section, I want to point to a few instances in Zimmer's own work where the conjecture was proven assuming slow enough growth of derivatives.  First in \cite{MR743815}, the conjecture is proven for ergodic volume preserving actions with what one might call {\em immediate subexponential growth of derivatives}, that is where we have
$$\|D\alpha(\gamma)\|\le e^{\epsilon \wl(\gamma)}.$$
\noindent This condition holds in particular for perturbations of isometric actions.  In the original proof, Zimmer
 used that $\Gamma$ was a higher rank lattice, but later he improved the result to cover all groups with property $(T)$ \cite{MR900826}.  A stronger result, not requiring either ergodic or volume preserving, can be proven using the techniques I developed with Margulis in \cite{MR2198325}.  In this essentially perturbative setting, my results with Margulis are clearly more robust than Zimmer's, allowing us to prove foliated results in \cite{MR2198325} with further applications in \cite{MR2521112}.  In a later paper, Zimmer introduced the notion of distal action and proved the conjecture for ergodic distal actions \cite{MR784532} by an argument similar to the one in \cite{MR743815}.  From the current point of view, being distal amounts to having a smooth invariant volume and having (eventual) polynomial growth of derivatives. For smooth volume preserving ergodic actions, Zimmer's proof here actually gives a proof assuming only subexponential growth of derivatives.  Another proof under those same hypotheses can be given by following Zimmer's arguments in \cite{MR1147291}. The proof above using strong property $(T)$ is both more robust by not requiring a smooth volume and considerably simpler than following Zimmer's arguments.  However, essentially all of Zimmer's to arguments in this context apply to ergodic volume preserving actions of groups with property $(T)$ which also preserve a measuable metric. Many hyperbolic groups, including lattices in $Sp(1,n)$ and random groups in many models, have property $(T)$ and fail to have strong property $(T)$ and Zimmer's results may be useful in the study of actions of those groups.  In any case it was already clear from Zimmer's earliest work that proving some kind of uniformly slow growth of derivatives was key to proving the conjecture.  What was missing until very recently was the ability to prove any estimate of that kind, an ability provided by the connection to hyperbolic dynamics, Lyapunov exponents, and rigidity of invariant measures in both the homogeneous and inhomogeneous setting.

As a last remark in this section, there is an unresolved issue in the current proofs that is obscured by the special case in which I state all results.  Namely, in the current versions of the argument by Brown, Rodriguez Hertz and Wang, the results one can prove only see the number of root spaces and not their multiplicities.  For this reason, for non-split groups, the current state of the art is often quite far from optimal.  My current sense is that while this difficulty is serious and non-trivial, it is one that will be overcome in time.

\section{Future directions}
\label{sec:future}

In this section, I will briefly discuss some of the most promising future directions and connections and end with
a discussion of some related work which I think may help point the way to a more general diffeomorphism group valued representation theory of finitely generated groups.

A key ingredient in the discussion above is close connections to homogeneous dynamics and the study of
invariant measures in that context. This area has been intensely  studied since the 60s or 70s with key contributions by Dani, Einseidler, Furstenberg, Katok, Lindenstrauss, Margulis, Raghunathan, Ratner, Shah, Spatzier  and Witte Morris.   Progress in that area has been quite dramatic since Ratner proved her theorem on the classification of measures invariant under unipotent flows in the early 90's and even our arguments in the case of non-uniform lattices barely scratch the surface. I want
to point to a certain thread of ideas particularly relevant to the recent work on Zimmer's conjecture. A key impetus
comes from work of Margulis and Tomanov where they, in the course of proving an extension of Ratner's theorems, introduce entropy techniques based in the work of Ledrappier and Young \cite{MR1253197,MR819557, MR819556, MR743818}.
The shorter argument in Cantat's paper for special cases of Theorem \ref{theorem:invariantmeasure} and \ref{theorem:invariantmetric}, uses ideas that appear in some form in those papers.  Another key development that began soon afterwards in work of Katok and Spatzier on invariant measures for higher rank abelian groups is to use a more refined form of the entropy argument, one that sees coarse Lyapunov subspaces and not just Lyapunov subspaces \cite{MR1406432}.  This idea is used constantly since that time in the study of invariant measures for higher rank abelian groups including the work of Einsiedler, Katok and Lindenstrauss \cite{MR2195133, MR2191228, MR2247967, MR3296819}.  For some summary of those ideas and a particularly accessible account of the entropy lemmas in the homogeneous setting, I recommend the survey by Einsiedler and Lindenstrauss \cite{MR2648695}.  A key contribution to the progress reported
in Section \ref{sec:zimmer} is the work of Brown, Rodriguez Hertz and Wang, which develops the entropy theory along coarse Lyapunov foliations in the general, non-homogeneous, setting \cite{AWB-GLY, AWBFRHZW-latticemeasure}.
 This can be seen as part of a promising set of  results that bring measure rigidity into the non-linear setting. The earliest of these results are for actions of abelian groups in the work of Katok, Kalinin and Rodriguez Hertz (see for example  \cite{MR2811602}) and these will likely prove relevant, but I want to place more emphasis on a different, more recent direction.

A major development in the homogeneous setting occured with work of Benoist and Quint \cite{BQIII}, where they study stationary and invariant measure for quite general groups.  These ideas then inspired the remarkable breakthrough work by Eskin and Mirzakhani on stationary and invariant measures for the $\Sl(2,\R)$ action on the Teichmuller moduli space \cite{1302.3320}, a setting that one might label semi-homogeneous. Using some ideas from \cite{1302.3320}, Brown and Rodriguez Hertz have written a magesterial paper about stationary and invariant measures for groups acting on surfaces subject to some natural hyperbolicity conditions \cite{MR3671937}. This is a surprising connection of rigidity of group actions to Teichmuller dynamics.   I hope the discussion of the arguments from \cite{BFH} convince readers that there is interest in hyperbolic measures even if one is not a priori interested in hyperbolic dynamics.  Roughly speaking their results show that if one has a hyperbolic stationary measure $\mu$ for a generating measure of some group $\Gamma$ acting by $C^2$ volume preserving diffeomorphisms on a compact surface $M$, one has a trichotomy

\begin{enumerate}
\item $\mu$ is supported on a finite set of points,
\item there is a $\Gamma$ invariant line field or
\item $\mu$ is, up to normalization, the restriction of volume to a $\Gamma$ invariant set.
\end{enumerate}

\noindent In particular, in the first and third case, $\mu$ is $\Gamma$ invariant. This result is potentially relevant to an another old conjecture of Zimmer's:  if $\Gamma$ is a group with property $(T)$ then any volume preserving action of $\Gamma$ on a surface factors through a finite quotient. In fact, if one can eliminate cae $(2)$ above for actions of groups with property $(T)$, then it is relatively easy to prove subexponential growth of derivatives.  While it is an easy consequence of a result of Zimmer that no measure preserving action of a group with property $(T)$ can preserve a line field on a surface, this isn't known for non-measure preserving actions, which makes application of the result somewhat difficult.

In joint work with Eskin, Brown and Rodriguez Hertz are developing a high dimensional analogue of this theorem which is somewhat more difficult to state.  This work should have further ramifications for work on rigidity of group actions in higher dimensions.

Brown, Rodriguez Hertz and Wang believe that the methods of \cite{AWBFRHZW-latticemeasure} can be adapted to show that if a lattice in $\Sl(n,\R)$ acts on a manifold of dimension $n-1$ without invariant measure then the action is smoothly conjugate to the standard one on a sphere $S^{n-1}$ or a projective space $P(\R^n)$.  Combining their ideas with the proof of Zimmer's conjectures in this case should completely classify all smooth actions of lattices in $\Sl(n,\R)$ on manifolds of dimension at most $n-1$.  The next test case for our techniques is clearly manifolds of dimension $n$ and while progress seems possible, this would be a major advance.  The key difficulty is that there are many more examples here. Not only can one take the $\Sl(n,\Z)$ action on the torus $\T^n$, but one can blow up periodic orbits for this action and also glue along blown up periodic orbits to produce manifolds of quite complex topology.  These examples were first discovered by Katok and Lewis and their properties were further explored in work I did with Benveniste and Whyte \cite{MR1380646, MR1866848, MR2154667, MR2342012}.  A reasonable first case to study is to assume that all stationary measures for the $\Gamma$ action are invariant and see if one can prove the action is smoothly conjugate to the standard action of $\Sl(n,\Z)$ on $\T^n$.  Another reasonable test case is to see if having a non-invariant stationary measure of full support implies the action is a skew product action on $P(\R^n)\times S^1$ or $S^{n-1} \times S^1$.  Many such actions can be constructed by viewing $P(\R^n)$ as $G/P$, taking a homomorphism $P\rightarrow \R$ and inducing any vector field on $S^1$ to an action of $G$ on $G/P \times S^1$.  It is also true that $\Sl(n, \R)$ acts on the $P(\R^{n+1})$ and $S^n$ and therefore so do all lattices in $\Sl(n,\R)$.  In this case there are orbits where the dynamics is dissipative instead of conservative and new ideas are likely needed.

I want to close this paper by discussing other developments that point much further into the future of the Zimmer program and its variants and generalizations.  Zimmer's program is essentially the generalization of the study of discrete subgroups of Lie groups and their finite dimensional representation theory to the study of those infinite dimensional representations that arise via actions on compact manifolds.  This analogy has proven very robust and Hurtado's Burnside Theorem is certainly a good example of another direction one can pursue.  Here I want to mention another, that is trying to study the {\em representation variety} or {\em character variety} in settings where there are many representations.  The word ``variety" should not be taken too seriously here as a principle difficulty of this setting is that there is no algebraic geometry that applies to groups of diffeomorphisms or homeomorphisms.  Despite this, in very recent work Mann and Wolff have given a complete characterization of the {\em rigid} representations or isolated points for the $\Homeo(S^1)$ character variety of the fundamental group of a compact surface \cite{MannWolff}.  These are exactly the representations that arise by restricting the action of a finite dimensional Lie subgroup of $\Homeo(S^1)$ to a cocompact lattice, a surprising and deep connection with the finite dimensional setting.  In this setting the finite dimensional Lie group that arises is necessarily a finite cover of $\PSl(2,\R)$.  This result is promising, but it is clearly a long way from understanding diffeomorphism or homeomorphism valued representation varieties in any detail.

  \bibliographystyle{AWBmath}

\bibliography{bibliography}

\end{document}